\documentclass[a4paper, english, 12pt, twoside]{article}
\usepackage{babel}
\usepackage[T1]{fontenc}
\usepackage{amssymb}
\usepackage{eufrak}
\usepackage[all]{xy}

\newcommand{\proof}{\noindent {\it Proof:}\ }

\begin{document}
\title{{\bf A Hopf algebra having a separable Galois extension is finite dimensional}}
\author{J. Cuadra\thanks{This research was supported by project
MTM2005-03227 from MCYT and FEDER.} \\
{\small Universidad de Almer\'{\i}a} \\
{\small Depto. \'{A}lgebra y An\'{a}lisis Matem\'atico} \\
{\small E-04120 Almer\'{\i}a, Spain} \\
{\small email: jcdiaz@ual.es}}
\date{ }
\maketitle


\begin{abstract}
It is shown that a Hopf algebra over a field admitting a Galois
extension separable over its subalgebra of coinvariants is of
finite dimension. This answers in the affirmative a question posed
by Beattie et al. in [{\it Proc. Amer. Math. Soc.} 128 No. 11
(2000), 3201-3203].
\end{abstract}

\section*{Introduction}

The notion of Hopf-Galois extension, as known nowadays, is due to
Kreimer and Takeuchi \cite{KT} and it is a mainstay of Hopf
algebra theory. It emanated from the work of Chase and Sweedler
about actions of Hopf algebras on rings \cite{CS}. When the Hopf
algebra is the coordinate algebra of an affine group scheme that
acts on an affine algebra, the above notion may be interpreted in
geometric terms and it is linked with the concept of torsor or
principal homogeneous space \cite[page 168]{Sch}. Faithfully flat
Hopf-Galois extensions are currently widely accepted as a
noncommutative counterpart of this geometric concept.
\par \smallskip

For $H=k[G]$, the group algebra of a group $G$ over a field $k$,
\cite[Theorem 8.1.7]{Mo} shows that an $H$-Galois extension is
precisely a strongly graded algebra. That is, a $k$-algebra $A$
admitting a decomposition $A=\oplus_{\sigma \in G} A_{\sigma}$ as
$k$-vector space and satisfying
$A_{\sigma}A_{\tau}=A_{\sigma\tau}$ for all $\sigma,\tau \in G$.
The subalgebra of coinvariants of $A$ is $A_e$ ($e$ the identity
element of $G$). N\u{a}st\u{a}sescu et al. characterized in
\cite[Proposition 2.1]{DVV} when the extension $A_e\subset A$ is
separable. In particular, they found that if $A$ is separable
over $A_e$, then $G$ is finite. It was investigated in \cite{BDR}
if an analogous result could hold for general Hopf algebras. To
be more precise, suppose that $H$ is a Hopf algebra over $k$
having a Hopf-Galois extension separable over its subalgebra of
coinvariants. Is $H$ necessarily of finite dimensional? A
positive response was given under the additional assumption that
$H$ to be co-Frobenius. \par \smallskip

In this short note we answer this question in the affirmative.
Our proof relies on a combination of the properties of the
separability idempotent, the Galois maps and an old result of
Sweedler. As often happens in Hopf algebra theory, the new proof
seems more natural and simpler than the original proof for
$H=k[G]$. Combining our result with one of Cohen and Fischman we
provide a characterization of separable Hopf-Galois extensions
that generalizes to Hopf algebras the above-mentioned one of
N\u{a}st\u{a}sescu et al. for strongly graded rings. \par
\bigskip

We fix some notation and recall the definition of Hopf-Galois
extension. We expect that the reader is familiar with the
rudiments of Hopf algebra theory. Our conventions and notations
are those of \cite{Mo}. Throughout $H$ stands for a Hopf algebra
over a field $k$. Its counit is denoted as usual by
$\varepsilon$. All vector spaces considered in the sequel are
over $k$, map means linear map, and $\otimes$ denotes the tensor
product over $k$. \par
\smallskip

For a right $H$-comodule algebra $A$ with structure map $\rho:A
\rightarrow A \otimes H$, its subalgebra of coinvariants
$A^{co(H)}=\{a \in A: \rho(a)=a \otimes 1_H\}$ is denoted by $B$.
The Galois maps are given by
\begin{eqnarray*}
can:A \otimes_B A \rightarrow A \otimes H, & & a \otimes_B a'
\mapsto \sum_{(a')} aa'_{(0)} \otimes a'_{(1)}, \\
can':A \otimes_B A \rightarrow A \otimes H, & & a \otimes_B a'
\mapsto \sum_{(a)} a_{(0)}a' \otimes a_{(1)}.
\end{eqnarray*}
Recall from \cite[Definition 8.1.1]{Mo} that $B \subset A$ is
said to be an $H$-Galois extension if $can$ is an isomorphism. It
is known that for $H$ having bijective antipode, any of the two
maps to be isomorphism may be required in the definition of Galois
extension since $can$ is bijective if and only if $can'$ is
bijective, \cite[page 124]{Mo}. However, we will not employ this
fact. Furthermore, {\it we will only use that $can$ is
surjective.}

\section{The main theorem}

\noindent {\bf Theorem.} {\it Let $H$ be a Hopf algebra and let
$A$ be a right $H$-Galois extension separable over its subalgebra
of coinvariants $B$. Then $H$ is finite dimensional.} \par
\medskip

\proof We will prove that $H$ has a non-zero finite dimensional
left ideal. In virtue of \cite[Corollary 2.7]{Sw2} this will imply
that $H$ is finite dimensional. \par \smallskip

Let $e=\sum_{i=1}^n e_i \otimes_B e'_i \in A \otimes_B A$ be the
separability idempotent given by hypothesis. Then $\sum_{i=1}^n
e_ie'_i=1_A$ and $ae=ea$ for all $a \in A.$ Notice that $can'(e)$
is non-zero since
$$1_A=\sum_{i=1}^n e_ie'_i=\sum_{i=1}^n \sum_{(e_i)} e_{i(0)}\varepsilon(e_{i(1)})
e'_i=(id_A \otimes \varepsilon)can'(e).$$ We pick non-zero elements
$a_j\in A, h_j \in H$ for $j=1,...,m$ such that
\begin{equation}\label{eq1}
\sum_{j=1}^m a_j \otimes h_j=can'(e)=\sum_{i=1}^n \sum_{(e_i)}
e_{i(0)}e'_i \otimes e_{i(1)}
\end{equation}
and the $a_j$'s are linearly independent. Take $h \in H$ arbitrary.
Since $can$ is surjective we may find $c_l, d_l \in A$ for
$l=1,...,r$ satisfying
\begin{equation}\label{eq2}
1 \otimes h=can(\sum_{l=1}^r c_l \otimes_B d_l)=\sum_{l=1}^r
\sum_{(d_l)} c_ld_{l(0)} \otimes d_{l(1)}.
\end{equation}
For $l=1,...,r$ we have $d_le=ed_l$. Applying $can'$ to this
equality we get
\begin{equation}\label{eq3}
\sum_{i=1}^n \sum_{(e_i)} \sum_{(d_l)} d_{l(0)}e_{i(0)}e'_i
\otimes d_{l(1)}e_{i(1)}=\sum_{i=1}^n \sum_{(e_i)}
e_{i(0)}e'_id_l \otimes e_{i(1)}.
\end{equation}
Then,
\begin{flushleft}
\begin{eqnarray*}
\sum_{j=1}^m a_j \otimes hh_j & \stackrel{(\ref{eq1})
(\ref{eq2})}{=} & \sum_{l=1}^r \sum_{(d_l)} \sum_{i=1}^n
\sum_{(e_i)} c_ld_{l(0)}e_{i(0)}e'_i \otimes
d_{l(1)}e_{i(1)} \smallskip \\
& \stackrel{(\ref{eq3})}{=} & \sum_{l=1}^r \sum_{i=1}^n
\sum_{(e_i)} c_le_{i(0)}e'_id_{l} \otimes e_{i(1)}  \smallskip \\
 & \stackrel{(\ref{eq1})}{=} & \sum_{l=1}^r \sum_{j=1}^m c_la_jd_l \otimes h_j.
\end{eqnarray*}
\end{flushleft}
Let $\varphi_t \in A^*$ be such that $\varphi_t(a_j)=\delta_{tj}$,
the Kronecker symbol, for $t=1,...,m.$ Evaluating $\varphi_t
\otimes id_H$ on the preceding set of equalities we obtain
$$hh_t=\sum_{l=1}^r \sum_{j=1}^m \varphi_t(c_la_jd_j) h_j.$$
This yields that the subspace spanned by the $h_j$'s is a finite
dimensional non-zero left ideal of $H$, as required. {\it q.e.d.}
\par \bigskip \smallskip

Cohen and Fischman provided in \cite[Theorem 1.8]{CF} several
characterizations of separable Hopf-Galois extensions for a
finite dimensional Hopf algebra. These characterizations together
with our result allow to characterize such extensions for an
arbitrary Hopf algebra.  \par \bigskip

\noindent {\bf Corollary.} {\it Let $H$ be a Hopf algebra and let
$A$ be a right $H$-Galois extension. Then, $A^{co(H)}\subset A$ is
separable if and only if $H$ is finite dimensional and one of the
equivalent conditions (2)-(6) in \cite[Theorem 1.8]{CF} holds.}
\par\bigskip

This corollary may be viewed as a generalization to Hopf algebras
of \cite[Proposition 2.1]{DVV} characterizing strongly graded
rings that are separable over its component of degree one.


\begin{thebibliography}{9}

\bibitem{BDR} M. Beattie, S. D\u{a}sc\u{a}lescu and \c{S}. Raianu, {\it A
Co-Frobenius Hopf Algebra with a Separable Galois Extension is
Finite}. Proc. Amer. Math. Soc. {\bf 128} No. 11 (2000),
3201-3203.

\bibitem{CS} S.U. Chase and M.E. Sweedler, {\it Hopf Algebras and
Galois Theory}. Lecture Notes in Mathematics {\bf 97}.
Springer-Verlag, Berlin, 1969.

\bibitem{CF} M. Cohen, D. Fischman, {\it Semisimple Extensions and
Elmenents of Trace 1}. J. Algebra {\bf 149} (1992), 419-437.

\bibitem{DVV} C. N\u{a}st\u{a}sescu, M. Van den Bergh and F. Van
Oystaeyen, {\it Separable Functors Applied to Graded Rings}. J.
Algebra {\bf 123} (1989), 397-413.

\bibitem{KT} H.F. Kreimer, M. Takeuchi, {\it Hopf Algebras and
Galois Extensions of an Algebra}. Indiana Univ. Math. J. {\bf 30}
(1981), 675-692.

\bibitem{Mo}  S. Montgomery,  {\it Hopf Algebras and Their Actions on
Rings}. CMBS No. 82, AMS, 1993.\vspace{-1ex}

\bibitem{Sch} H.-J. Schneider, {\it Principal Homogeneous Spaces
for Arbitrary Hopf algebras}. Israel J. Math. {\bf 72} Nos. 1-2
(1980), 167-195.

\bibitem{Sw2} M.E. Sweedler, {\it Integrals for Hopf Algebras}.
Ann. of Math. (2) {\bf 89} (1969), 323-335.
\end{thebibliography}
\end{document}